\documentclass[12pt]{amsart}
\usepackage{amsmath}
\usepackage{amsfonts}
\usepackage{amssymb}
\usepackage[all]{xy}

\newcommand{\ZZ}{\mathbb{Z}}
\newcommand{\FF}{\mathbb{F}}
\newcommand{\CC}{\mathbb{C}}
\newcommand{\QQ}{\mathbb{Q}}

\newcommand{\KK}{\mathbb{K}}
\newcommand{\GG}{\mathbb{G}}
\newcommand{\A}{\mathbb{A}}

\newcommand{\Gm}{\GG_m}
\newcommand{\Ga}{\GG_a}
\newcommand{\Cm}{C^{\otimes m}}
\newcommand{\Cn}{C^{\otimes n}}

\DeclareMathOperator{\Hom}{Hom}
\DeclareMathOperator{\Der}{Der}
\DeclareMathOperator{\DR}{DR}
\DeclareMathOperator{\Ext}{Ext}
\DeclareMathOperator{\Exp}{Exp}
\DeclareMathOperator{\Lie}{Lie}
\DeclareMathOperator{\Mat}{Mat}
\DeclareMathOperator{\wt}{wt}

\DeclareMathOperator{\Spec}{Spec}

\newcommand{\assign}{\mathrel{:=}}
\newcommand{\laurent}[2]{#1(\!( #2 )\!)}

\newcommand{\smfrac}[2]{{\textstyle \frac{#1}{#2}}}
\newcommand{\tpi}{\widetilde{\pi}}
\newcommand{\iso}{\stackrel{\sim}{\rightarrow}}

\theoremstyle{plain}
\newtheorem{theorem}{Theorem}[section]
\newtheorem{lemma}[theorem]{Lemma}
\newtheorem{corollary}[theorem]{Corollary}
\newtheorem{proposition}[theorem]{Proposition}

\theoremstyle{definition}

\theoremstyle{remark}
\newtheorem{example}[theorem]{Example}
\newtheorem{remark}[theorem]{Remark}
\newtheorem*{notation}{Notation}

\begin{document}


  \title{A Weil-Barsotti formula for Drinfeld modules}
  \author{Matthew A.~Papanikolas}
  \address{Department of Mathematics\\
        Brown University\\
        Providence, RI 02912 USA}
  \email{map@math.brown.edu}
  \author{Niranjan Ramachandran}
  \address{Department of Mathematics\\
    University of Maryland\\
    College Park, MD 20742 USA}
  \address{Max-Planck-Institut f\"ur Mathematik\\
    Vivatsgasse 7\\ D-53111 Bonn, Germany}
  \email{atma@math.umd.edu}

  \date{April 17, 2002}

\begin{abstract}
  We study the group of extensions in the category of Drinfeld modules
  and Anderson's $t$-modules, and we show in certain cases that this
  group can itself be given the structure of a $t$-module.  Our main
  result is a Drinfeld module analogue of the Weil-Barsotti formula
  for abelian varieties.  Extensions of general $t$-modules are also
  considered, in particular extensions of tensor powers of the Carlitz
  module.  We motivate these results from various directions and
  compare to the situation of elliptic curves.
\end{abstract}


\maketitle


\section{Introduction and statement of results} \label{S:intro}

In this paper, we investigate extensions of Drinfeld modules using the
well-known analogy between abelian varieties and Drinfeld modules. We
prove analogues for Drinfeld modules of the classical Weil-Barsotti
formula and the Cartier-Nishi biduality theorem for abelian varieties.

Let $A$ be an abelian variety over a field $k$; we denote the dual
abelian variety by $A^{\vee}$. The Weil-Barsotti formula states that
for any $k$-algebra $R$, there is a natural, functorial isomorphism
$\Ext^1_R(A,\Gm) \cong A^{\vee}(R)$ where the first group is
calculated in the category of group schemes over $\Spec R$
\cite{milne86AV}. In other words, the functor $R \mapsto
\Ext^1_R(A,\Gm)$ (on $k$-algebras) is represented by the dual abelian
variety $A^{\vee}$.  The \emph{biduality} theorem of Cartier-Nishi
states that there is a canonical isomorphism $\Ext^1_R(A^{\vee},\Gm)
\cong A(R)$; this can also be restated: there is a canonical
isomorphism of abelian varieties $(A^{\vee})^{\vee} \cong A$.
Theorem~\ref{T:WB} below provides a Drinfeld module analogue of these
results.

There is also an important relationship between the de Rham cohomology
(in characteristic zero) of $A$ and the universal additive (or
vectorial) extension $A^{\natural}$ of $A$ \cite{mes}.  A de Rham
theory for Drinfeld modules based on additive extensions has already
been developed by G.~Anderson, P.~Deligne, E.-U.~Gekeler, and J.~Yu
\cite{gek89}. We pursue generalizations for $t$-modules.

We remark that extensions by $\Gm$ and $\Ga$ have been used by Deligne
in \cite[\S 10]{deli} to define Cartier duality and the de Rham theory
for $1$-motives.

\begin{notation}
  Let $K$ be a perfect field of characteristic $p > 0$, and let
  $\FF_q[t]$ be the polynomial ring in one variable over the finite
  field $\FF_q$ where $q = p^m$.  Fix an $\FF_q$-linear homomorphism
  $\iota\colon \FF_q[t] \to K$ with $\theta \assign \iota(t)$.
  Throughout, all Drinfeld modules and $t$-modules are defined with
  respect to the map $\iota$, and in particular all Drinfeld modules
  are $\FF_q[t]$-modules.
  
  The ring $K\{\tau\}$ is the ring of twisted polynomials in $K$ such
  that for $x \in K$, $\tau x = x^q \tau$.  A $d$-dimensional
  $t$-module over $K$ is at first an $\FF_q$-linear ring homomorphism
\[
  \Phi : \FF_q[t] \to \Mat_d(K\{\tau\}),
\]
such that, as a polynomial in $\tau$ with coefficients in $\Mat_d(K)$,
\[
  \Phi(t) = (\theta I_d + N)\tau^0 + M_1 \tau^1 + \cdots,
\]
where $I_d$ is the identity matrix and $N$ is nilpotent.  In general,
a $t$-module over $K$ is an algebraic group $E$ defined over $K$,
which is isomorphic over $K$ to $\Ga^d$, together with a choice of
$\FF_q$-linear endomorphism $t\colon E \to E$ such that $d(t-\theta)^n
\Lie(E) = 0$ for all $n$ sufficiently large.  By choosing an
isomorphism $E \cong \Ga^d$, one can specify a homomorphism $\Phi\colon
\FF_q[t] \to \Mat_d(K\{\tau\})$ as above.  To denote this choice of
coordinates, we write $E = (\Ga^d,\Phi)$.

Let $C$ denote the Carlitz module, $C \colon \FF_q[t] \to K\{\tau\}$,
defined by $C(t) = \theta + \tau$.

We take $\Ext^1(\cdot,\cdot)$ to be the bifunctor $\Ext^1$ from the
additive category of $t$-modules to the category of abelian groups.
In \S\ref{S:ExtBid} we see that, for two $t$-modules $E$ and $F$,
those extensions which induce trivial $t$-module extensions of their
respective tangent spaces comprise a canonical subgroup $\Ext^1_0(E,F)
\subseteq \Ext^1(E,F)$.  For a $t$-module $E$, we let $E^{\vee}
\assign \Ext^1_0(E,C)$.
\end{notation}

Our analogue of the classical Weil-Barsotti formula and the
Cartier-Nishi biduality theorem is the following  
\begin{theorem} \label{T:WB}
Let $E$ be a Drinfeld module of rank $r \geq 2$.
\begin{enumerate}
\item[(a)] The group $\Ext^1(E,C)$ is naturally a $t$-module of
  dimension $r$ and sits in an exact sequence of $t$-modules
\[
  0 \to E^{\vee} \to \Ext^1(E,C) \to \Ga \to 0.
\]
Furthermore, $E^{\vee}$ is the Cartier-Taguchi dual $t$-module
associated to $E$ \cite{tag95}, and in particular, $E^{\vee}$ is
isomorphic to the $(r-1)$-st exterior power $\wedge^{r-1}E$ of
$E$.
\item[(b)] The group $\Ext^1(E^{\vee},C)$ is also naturally a
  $t$-module of dimension $r$ and sits in an exact sequence
\[
  0 \to E \to \Ext^1(E^{\vee},C) \to \Ga^{r-1} \to 0.
\]
Moreover, we have a biduality: $(E^{\vee})^{\vee} \cong E$.

\item[(c)] Any morphism $\beta:E \to F$ of Drinfeld modules (of
  rank $\geq 2$) induces a morphism of dual $t$-modules
  $\beta^{\vee}: F^{\vee} \to E^{\vee}$.
\end{enumerate}
\end{theorem}

The proof of Theorem~\ref{T:WB} also shows that the $t$-module
structure on $\Ext^1(E,C)$ is compatible with base change of the field
$K$; see \S\ref{S:PerExt}.

Parts (a) and (b) of Theorem~\ref{T:WB} for Drinfeld modules of rank
$2$ have been proven by S.~S.~Woo \cite{woo95}.  Y.~Taguchi
\cite{tag95} has constructed a Weil pairing (compatible with the
Galois action) on the torsion points of $E$ and $E^{\vee}$.  Taguchi
remarks in \cite{tag95} that his definition of the Cartier dual
$E^{\vee}$ of a Drinfeld module $E$ does generalize to some (but not
all) $t$-modules.

Theorem~\ref{T:WB} requires us to work outside the category of Drinfeld
modules, and one may ask for general $t$-modules $E$ and $F$ over $K$
whether $\Ext^1(E,F)$ has the structure of a $t$-module.  In this vein
we have the following result.  Let $\Cn$ denote the $n$-th tensor
power of the Carlitz module \cite{andthk90}.

\begin{theorem} \label{T:Ctens}
  If $n > m$, then $\Ext^1(\Cm,\Cn)$ has the structure of a
  $t$-module, and there is an exact sequence of $t$-modules
\[
  0 \to C^{\otimes(n-m)} \to \Ext^1(\Cm,\Cn) \to L \to 0,
\]
where $L$ is an $m$-dimensional iterated extension of $\Ga$.
Moreover, $\Ext^1_0(\Cm,\Cn) \cong C^{\otimes(n-m)}$.
\end{theorem}

Since the (tractable) period $\tpi^{n-m}$ of $C^{\otimes(n-m)}$ is a
power of the period $\tpi$ of the Carlitz module (see D.~Goss
\cite[Ch.~3]{goss:FF}), one should compare Theorem~\ref{T:Ctens} with
the isomorphism
\[
  \Ext^1_{\mathrm{MHS}}(\ZZ(m),\ZZ(n)) \cong \CC/(2\pi i)^{n-m} \ZZ,
  \quad n > m,
\]
{from} the theory of mixed Hodge structures \cite{carl}.

The question of whether an analogue of the Weil-Barsotti formula holds
for general $t$-modules is also interesting.  Experimental evidence
suggests that the general situation is subtle and that such formulas
are not always valid for pure $t$-modules, e.g.\ in the form of
Theorem~\ref{T:WB}, a Weil-Barsotti formula does not hold for $\Cm$
because $\Ext^1(\Cm,\Cn)$ is not well-behaved for $n \leq m$.
See \S\ref{S:Ctens} for more details.

As pointed out by the referee, it would be worth investigating the
extent to which Theorems~\ref{T:WB} and~\ref{T:Ctens} are true for
Drinfeld modules over rings more general than $\FF_q[t]$.  This raises
some technical issues, which we discuss in \S\ref{S:PerExt}.

The outline of this paper is as follows.  In \S\ref{S:ExtBid}, we
present definitions and fundamental results on extensions of
$t$-modules.  We prove Theorems~\ref{T:WB} and~\ref{T:Ctens} in
\S\S\ref{S:WB}--\ref{S:Ctens}.  In \S\ref{S:PerExt}, we consider
extensions of $t$-modules from an analytic viewpoint, so as to
motivate the expectation that $\Ext^1(E,F)$ can be represented by a
$t$-module for certain $t$-modules $E$ and $F$.  We consider the
situation of elliptic curves in \S\ref{licht} and compare our results
to an unpublished theorem of S.~Lichtenbaum about extensions
of elliptic curves over $\CC$.  We conclude in \S\ref{S:Rmk} with some
remarks about extensions of $t$-motives.

\section{Extensions of $t$-modules and biderivations} \label{S:ExtBid}

In this section we establish definitions and results about extensions
of $t$-modules.  For general definitions of $t$-modules we follow the
terminology in Goss \cite[Ch.~5]{goss:FF}.  So as not to lead to
confusion, we adhere to the following convention: a ``$t$-module''
refers to the object of the same name defined in \S\ref{S:intro},
whereas an ``$\FF_q[t]$-module'' is simply a module over the ring
$\FF_q[t]$.

We point out that the results in this section remain valid in the case
that $K$ is not perfect, though we do not make use of this fact
later on.

Let $E$ and $F$ be $t$-modules over $K$.  An extension of $E$ by $F$
is a $t$-module $X$ fitting into an exact sequence of $t$-modules
\begin{equation} \label{E:FXE}
  0 \to F \to X \to E \to 0.
\end{equation}
Then $\Ext^1(E,F)$ is defined to be the group (under Baer sum) of of
$t$-module extensions of $E$ by $F$ up to Yoneda equivalence.  

The main tool which enables us to compute this group is that of
biderivations.  The following definitions run parallel to those of
W.~D.\ Brownawell and the first author \cite{bp:Gmm} and Gekeler
\cite{gek89}, where extensions of $t$-modules by $\Ga$ were
investigated.

Let $\Phi\colon \FF_q[t] \to \Mat_d(K)\{ \tau \}$ and $\Psi\colon
\FF_q[t] \to \Mat_e(K)\{ \tau \}$ be choices of coordinates for $E$
and $F$ respectively, where $\Mat_{d}(K)\{\tau\}$ is the ring of
twisted polynomials with matrix coefficients.  A
\emph{$(\Phi,\Psi)$-biderivation} is an $\FF_q$-linear map
\[
  \delta : \FF_q[t] \to \Mat_{e \times d}(K)\{ \tau \},
\]
which satisfies
\begin{equation} \label{E:deltab}
  \delta(ab) = \Psi(a)\delta(b) + \delta(a)\Phi(b), \quad
  \forall a, b \in \FF_q[t].
\end{equation}
The $\FF_q$-vector space of all $(\Phi,\Psi)$-biderivations is denoted
$\Der(\Phi,\Psi)$.  It is straightforward to check that a biderivation
$\delta$ is uniquely determined by the single value $\delta(t)$, and
so if $V \in \Mat_{e \times d}(K) \{ \tau \}$, we define $\delta_V \in
\Der(\Phi,\Psi)$ to be that biderivation such that $\delta_V(t) = V$.
In this way, we have an isomorphism of $\FF_q$-vector spaces
\begin{equation} \label{E:MatDer}
  V \mapsto \delta_V : \Mat_{e \times d}(K) \{ \tau \} \iso
  \Der(\Phi,\Psi).
\end{equation}
A biderivation $\delta$ is called \emph{inner} if for some $U \in
\Mat_{e \times d}(K\{ \tau \})$ we have
\begin{equation} \label{E:delin}
  \delta(a) = \delta^{(U)}(a) \assign U\Phi(a) - \Psi(a)U, \quad
  \forall a \in \FF_q[t].
\end{equation}
The subspace of $\Der(\Phi,\Psi)$ of inner biderivations is denoted
$\Der_{in}(\Phi,\Psi)$.

Every $(\Phi,\Psi)$-biderivation $\delta$ gives rise to an extension
$X = (\Ga^{d+e},\Upsilon)$ of $E$ by $F$ by defining
\[
  \Upsilon(a) \assign \begin{pmatrix}
    \Phi(a) & 0 \\
    \delta(a) & \Psi(a)
    \end{pmatrix}, \quad \forall a \in \FF_q[t].
\]
Again it is straightforward, using \eqref{E:deltab}, to check that
$\Upsilon$ is well-defined.  Moreover, every extension of $E$ by $F$
defines a unique biderivation.

We note that if $\delta^{(U)}$ is an inner biderivation then in fact $X$ is
split.  In this case the matrix $\Theta \assign \left( \begin{smallmatrix}
    I_d & 0 \\ U & I_e
  \end{smallmatrix} \right)$ provides the splitting, where $I_d$,
$I_e$ are identity matrices:
\[
\Theta^{-1} \Upsilon(a) \Theta = \begin{pmatrix}
  \Phi(a) & 0 \\
  0 & \Psi(a)
    \end{pmatrix}, \quad \forall a \in \FF_q[t].
\]
Furthermore, it follows from the above discussion that every split
extension arises in this way.  

Suppose we are given two extensions of $E$ by $F$ which are Yoneda
equivalent. It follows easily from the definition of Yoneda
equivalence that the corresponding biderivations differ by an inner
biderivation.  It is straightforward to check that the (Baer) sum on
$\Ext^1(E,F)$ corresponds to usual addition on the level of
biderivations.

Now the endomorphisms of $E$ and $F$ induce (identical)
$\FF_q[t]$-module structures on $\Ext^1(E,F)$.  That is, if $X$
represents a class in $\Ext^1(E,F)$ and $b \in \FF_q[t]$, we can
define two $t$-modules $X \ast b$ and $b \ast X$, which ultimately
represent the same class in $\Ext^1(E,F)$.  Explicitly, suppose
$\delta$ is the $(\Phi,\Psi)$-biderivation corresponding to $X$ and
$\pi\colon X \to E$ is the natural map in \eqref{E:FXE}.  Let
\[
  X \ast b \assign \ker \left( (e,x) \mapsto \Phi(b)(e) - \pi(x) : E
  \oplus X \to E \right).
\]
Then $X \ast b$ is itself a $t$-module extension of $E$ by $F$, and
the operation of $\FF_q[t]$ on it is given by
\[
  (\Upsilon \ast b)(a) = \begin{pmatrix}
    \Phi(a) & 0 \\
    \delta(a)\Phi(b) & \Psi(a)
    \end{pmatrix}, \quad \forall a \in \FF_q[t].
\]
On the other hand, we can similarly use endomorphisms of $F$ to define
an extension $b\ast X$ whose $t$-module structure is given by
\[
  (b \ast \Upsilon)(a) = \begin{pmatrix}
    \Phi(a) & 0 \\
    \Psi(b)\delta(a) & \Psi(a)
    \end{pmatrix}, \quad \forall a \in \FF_q[t].
\]
To see that $X\ast b$ and $b \ast X$ are equivalent extensions, we
note that $\varepsilon \colon \FF_q[t] \to \Mat_{e \times d}(K)\{ \tau
\}$ defined by
\[
  \varepsilon(a) \assign \delta(a)\Phi(b) - \Psi(b)\delta(a), \quad
  \forall a \in \FF_q[t],
\]
is in fact the inner biderivation $\delta^{(U)}$, with $U = \delta(t)$
in \eqref{E:delin}.  That is,
\begin{equation} \label{E:tDer}
  b \colon \delta(\cdot) \mapsto \delta(\cdot)\Phi(b)
  \quad \textnormal{and} \quad
  b \colon \delta(\cdot) \mapsto \Psi(b)\delta(\cdot)
\end{equation}
define $\FF_q[t]$-module structures on $\Der(\Phi,\Psi)$ which are the
same modulo $\Der_{in}(\Phi,\Psi)$.  We record the results from the
preceding paragraphs in the following lemma.

\begin{lemma} \label{L:ExtDer}
Let $E = (\Ga^d,\Phi)$ and $F = (\Ga^e,\Psi)$ be $t$-modules.  Then
\[
  \Ext^1(E,F) \cong \Der(\Phi,\Psi)/\Der_{in}(\Phi,\Psi)
\]
as $\FF_q[t]$-modules.
\end{lemma}

For $U \in \Mat_{d_1 \times d_2}(K)\{ \tau\}$, we let $dU \in
\Mat_{d_1 \times d_2}(K)$ be the constant term of $U$ as a
polynomial in $\tau$, and we define the following subspaces of
$\Der(\Phi,\Psi)$:
\begin{gather*}
  \Der_0(\Phi,\Psi) \assign \{ \delta \in \Der(\Phi,\Psi) : d\delta(t)
  = 0 \}, \\
  \Der_{si}(\Phi,\Psi) \assign \{ \delta^{(U)} \in \Der_{in}(\Phi,\Psi) :
  dU = 0 \}.
\end{gather*}
The utility of $\Der_0(\Phi,\Psi)$ is derived from the following
lemma, whose immediate corollary follows from Lemma~\ref{L:ExtDer}.
Note that $\Der_0$ represents a different object here than in
\cite{bp:Gmm} and \cite{gek89}.  Biderivations in
$\Der_{si}(\Phi,\Psi)$ are called \emph{strictly inner}, and clearly
$\Der_{si}(\Phi,\Psi) \subseteq \Der_0(\Phi, \Psi)$.  We will study
$\Der_{si}(\Phi,\Psi)$ in more detail in \S\ref{S:PerExt}.

The map $d\Phi \colon \FF_q[t] \to \Mat_d(K)$ defines a non-abelian
$t$-module whose underlying space is the tangent space $\Lie(E) \cong
K^d$.  Furthermore, the map
\[
  \delta \mapsto d\delta :
  \Der(\Phi,\Psi) \to \Der(d\Phi,d\Psi)
\]
is $\FF_q$-linear, and it is $\FF_q[t]$-linear modulo inner
biderivations.

\begin{lemma} \label{L:Dseq}
  Let $E = (\Ga^d,\Phi)$ and $F = (\Ga^e,\Psi)$ be $t$-modules.  The
  following is an exact sequence of $\FF_q[t]$-modules:
\begin{equation} \label{E:Dseq}
  0 \to \frac{\Der_0(\Phi,\Psi)}{\Der_0(\Phi,\Psi) \cap
  \Der_{in}(\Phi,\Psi)}
  \to \frac{\Der(\Phi,\Psi)}{\Der_{in}(\Phi,\Psi)}
  \to \frac{\Der(d\Phi,d\Psi)}{\Der_{in}(d\Phi,d\Psi)}.
\end{equation}
If $\theta \in K$ is transcendental over $\FF_q$, then the final map
in this sequence is surjective.
\end{lemma}

We define $\Ext^1_0(E,F)$ to be the $\FF_q[t]$-submodule of
$\Ext^1(E,F)$ corresponding to $\Der_0(\Phi,\Psi)/(\Der_0(\Phi,\Psi) \cap
  \Der_{in}(\Phi,\Psi))$.

\begin{corollary} \label{C:Eseq}
The sequence in Lemma~\ref{L:Dseq} corresponds to an exact sequence of
$\FF_q[t]$-modules,
\[
0 \to \Ext^1_0(E,F) \to \Ext^1(E,F) \to \Ext^1(\Lie(E),\Lie(F)),
\]
where the final map is surjective if $\theta \in K$ is transcendental
over $\FF_q$.
\end{corollary}

\begin{proof}[Proof of Lemma~\ref{L:Dseq}]
  Injectivity on the left of \eqref{E:Dseq} is clear.  To show
  exactness in the center, first any $\delta \in \Der_0(\Phi,\Psi)$
  maps to $0$ in $\Der(d\Phi,d\Psi)$.  On the other hand, suppose
  $\delta \in \Der(\Phi,\Psi)$ and $d\delta \in \Der(d\Phi,d\Psi)$ is
  inner, say $d\delta(t) = Ud\Phi(t) - d\Psi(t)U$ with $U \in
  \Mat_{e\times d}(K)\{ \tau \}$.  Then $\delta -
  \delta^{(U)}_{\Phi,\Psi} \in \Der_0(\Phi,\Psi)$ represents the same
  class as $\delta$ in $\Der(\Phi,\Psi)/\Der_{in}(\Phi,\Psi)$.
  
  In the case that $\theta \in K$ is transcendental over $\FF_q$, we
  show surjectivity on the right.  We suppose $d\Phi(t) = I_d \theta +
  M$ and $d\Psi(t) = I_e \theta + N$, where $M$ and $N$ are nilpotent.
  Without loss of generality we can assume $M$ and $N$ are both upper
  triangular.  The spaces of biderivations $\Der(d\Phi,d\Psi)$ and
  $\Der_{in}(d\Phi,d\Psi)$ are both naturally $K$-linear, and the map
  defining elements of $\Der_{in}(d\Phi,d\Psi)$,
\begin{equation} \label{E:UdelU}
  U \mapsto \delta^{(U)}_{d\Phi,d\Psi}(t) : \Mat_{e \times d}(K)\{\tau\} \to
  \Mat_{e \times d}(K)\{\tau\},
\end{equation}
is $K$-linear and respects the grading by degrees in $\tau$.  Now
\begin{equation} \label{E:ddelU}
  \delta^{(U)}_{d\Phi,d\Psi}(t) = U\theta - \theta U + UM - NU.
\end{equation}
If $U_{ij}$ is the matrix in $\Mat_{e \times d}(K)$ with a $1$ in the
$ij$-th entry and zeros elsewhere, then $\{ U_{m1}\tau^r,
U_{m-1,1}\tau^r, \ldots \} = \{ U_{ij} \tau^r \}_{i=m, \dots, e}^{j=1,
  \dots, d}$ is an ordered $K$-basis for $\Mat_{e \times d}(K)\tau^r$.
Using this basis, it follows from \eqref{E:ddelU} and the fact that
$M$ and $N$ are both upper triangular that the map in \eqref{E:UdelU},
restricted to $\Mat_{e \times d}(K) \tau^r$, is lower triangular with
$\theta^{q^r} - \theta$ along the diagonal.  Since $\theta^{q^r} -
\theta$ is non-zero by our assumption on $\theta$,
\[
  U \mapsto \delta^{(U)}_{d\Phi,d\Psi}(t) : \Mat_{e\times d}(K)\{ \tau
  \} \tau \iso \Mat_{e \times d}(K)\{ \tau \} \tau,
\]
is an isomorphism of $K$-vector spaces.  Therefore, for every $\delta
\in \Der(d\Phi,d\Psi)$, if $d\delta(t) = 0$, then $\delta$ is inner,
and so the right-hand side of \eqref{E:Dseq} is surjective.
\end{proof}

\section{A Weil-Barsotti formula} \label{S:WB}

In this section, we prove Theorem~\ref{T:WB}.  Let $E = (\Ga,\Phi)$ be
a Drinfeld module of rank $r \ge 2$, where $\Phi(t) = \theta + a_1
\tau + \dots + a_r \tau^r$.  For an extension $X$ of $E$ by
$C$, we let $\delta_X \in \Der(\Phi,C)$ be its associated biderivation.

\begin{proof}[Proof of Theorem~\ref{T:WB}(a)]
  An extension $X$ of $E$ by $C$ splits if and only if $\delta_X \in
  \Der_{in}(\Phi,C)$, i.e.\ if there exists $u \in K \{ \tau \}$ such
  that $\delta_X(t) = u \Phi(t) - (\theta + \tau) u$.  If $u = c
  \tau^m$, then we obtain an inner biderivation $\delta^{(u)}$ such that
\begin{equation} \label{bonn}
\begin{aligned}
  \delta^{(u)}(t) &= u \Phi(t) - (\theta + \tau) u \\
  &= c \theta^{q^m} {\tau}^m + c {a_1}^{q^m} {\tau}^{m+1} + \dots + c
  {a_r}^{q^m} {\tau}^{m+r} - \theta c{\tau}^m - c^q{\tau}^{m+1}.
\end{aligned}
\end{equation}
If $n \assign \deg_{\tau}(\delta_X(t))$ is greater than $r-1$, we can
repeatedly subtract (\ref{bonn}) from $\delta_X(t)$, with $m = n-r$,
$n-r-1, \dots, 0$, to reduce the $\tau$-degree of $\delta_X(t)$;
eventually this degree will be $<r$.  Namely, any
extension $X$ is equivalent to an extension $X'$ with
$\deg_{\tau}(\delta_{X'}(t)) \leq r-1$, which we call the reduced
representative of $X$.

According to \eqref{bonn}, the non-zero inner biderivations of least
degree have degree $r$, so two extensions $Z$ and $Y$ satisfying
$\deg_{\tau}(\delta_{Z}(t)) \leq r-1$, $\deg_{\tau}(\delta_{Y}(t))
\leq r-1$, and $\delta_{Z}(t) \neq \delta_{Y}(t)$, are inequivalent.
Therefore the map
\begin{equation} \label{E:Xred}
X \mapsto \delta_{X'}(t) :  \Ext^1(E,C) \iso V \assign \{ b_0 + b_1
  \tau + \dots + b_{r-1} \tau^{r-1} : b_i \in K\},
\end{equation}
where $X'$ is the reduced representative of $X$, induces an
isomorphism of $\FF_q[t]$-modules.

We now turn to the $t$-module structure on $\Ext^1(E,C)$.  Recall from
\eqref{E:tDer} that multiplication-by-$t$ on $\Ext^1(E,C)$ is defined
by $t\ast \alpha = (\theta + \tau) \alpha = \alpha\Phi(t)$.  Here we think
of $\alpha$ as being an element of $V$.  In order to see the action
explicitly, it is enough to consider $\alpha = b_i \tau^i$:
\begin{equation} \label{bson}
\begin{aligned}
t\ast (b_0) &= (\theta + \tau) b_0 = \theta b_0 + b_0^q \tau,\\
t\ast (b_1 \tau) & =   (\theta + \tau) b_1 \tau = \theta b_1 \tau + b_1^q
  \tau^2,\\ 
\vdots \quad & \phantom{=}\qquad \vdots\\
t\ast (b_{r-1} \tau^{r-1}) & =  (\theta + \tau)  b_{r-1} \tau^{r-1} =
  \theta b_{r-1} \tau^{r-1} + b_{r-1}^q \tau^r.
\end{aligned}
\end{equation}  
Using \eqref{bonn} with $m=0$, $c = {b_{r-1}^q}/{a_r}$, we can rewrite
the last identity of \eqref{bson} as
\begin{multline*}
  t\ast (b_{r-1} \tau^{r-1}) = \left(\smfrac{b_{r-1}^{q^2}}{a_r^q} -
    \smfrac{b_{r-1}^q a_1}{a_r}\right) \tau
  -\smfrac{b_{r-1}^q}{a_r}\left(a_2 \tau^2
    + \dots + a_{r-2} \tau^{r-2}\right) \\
  + \left(\theta b_{r-1} - \smfrac{b_{r-1}^q a_{r-1}}{a_r}\right)
  \tau^{r-1}.
\end{multline*}
Thus in terms of the elements $e_i = \tau^i$ ($i = 0, \dots, r-1$) of
$V$, the $t$-module structure on $\Ext^1(E,C)$ can be expressed by the
map $\Pi\colon \FF_q[t] \to \Mat_r(K)\{ \tau \}$, defined by
\[
  \Pi(t) \assign \begin{pmatrix}
    \theta & 0 & \cdots & \cdots & 0 & 0 \\
    \tau & \theta & 0 &  & \vdots  & -\frac{a_1}{a_r} \tau +
    \frac{1}{a_r^q}\tau^2 \\
    0  & \tau  & \theta & 0 & \vdots & -\frac{a_2}{a_r}\tau \\
\vdots & 0 & \tau & \theta  & 0  & \vdots\\
\vdots & & 0 & \tau & \theta & -\frac{a_{r-2}}{a_r}\tau\\
0 & \cdots &  \cdots & 0 & \tau & \theta -\frac{a_{r-1}}{a_r}\tau.
  \end{pmatrix}.
\]
Comparing with Taguchi in \cite[\S 5]{tag95}, it is clear that the
$t$-module $\Ext^1(E,C)$ is an extension of $\Ga$ by the $t$-module
denoted there $\check{E}$.  By Corollary~\ref{C:Eseq} and the
characterization of $V$ in \eqref{E:Xred}, it is clear $\Ext^1(E,C)$
is an extension of $\Ga$ by $E^{\vee}$.  Thus $E^{\vee}$ is the same
as Taguchi's $t$-module.  Moreover, Taguchi shows
\cite[Thm.~5.1]{tag95} that $E^{\vee}$ is a pure $t$-module isomorphic
to $\wedge^{r-1}E$.
\end{proof}

We turn to part (b), and for simplicity we assume that $a_r =1$; the
general case follows similarly.  The $t$-module structure on
$E^{\vee}$ is then defined by
\[
  \Psi(t) \assign \begin{pmatrix}
 \theta & 0  & \cdots & 0 & -{a_1}\tau + \tau^2   \\
 \tau  & \theta & 0 & \vdots & -{a_2}\tau \\
 0 & \tau & \theta & 0 & \vdots\\
\vdots & 0  & \tau & \theta & -{a_{r-2}}\tau\\
0 &  \cdots & 0 & \tau & \theta -{a_{r-1}}\tau.
  \end{pmatrix}.
\]

\begin{proof}[Proof of Theorem~\ref{T:WB}(b)]
  Consider the biderivation $\delta_X \colon \FF_q[t] \to \Mat_{1
    \times (r-1)}(K)\{ \tau \}$ of any extension $X$ of $E^{\vee}$ by
  $C$.  The inner biderivations are of the form $\delta^{(U)}(t) = U \Psi(t)
  - (\theta + \tau) U$ for $U = (u_i) \in \Mat_{1 \times
    ({r-1})}(K)\{\tau\}$.  Explicitly, $v = (v_i)
  \in \Mat_{1 \times (r-1)}(K)\{\tau \}$ defines an inner
  biderivation, $\delta^{(U)} = \delta_v$, if it has the form
\[
\begin{aligned}
v_i &= u_i \theta - \theta u_i + u_{i+1} \tau - \tau u_i \quad (1
\le i \le r-2), \\
v_{r-1} &= u_{r-1} \theta - \theta u_{r-1} + u_1 \tau^2 - \tau
u_{r-1} - {\textstyle \sum_{j=1}^{r-1} u_j a_j \tau}.
\end{aligned}
\]
An inner biderivation $\delta_v$ is said to be basic if $u_i = 0$ for all $i
\neq s$ and $u_s = c\tau^m$ with $c \in K$; we write $v = v(s,c,m)=
(v_1, \dots, v_{r-1})$. Explicitly written, these are 
\begin{equation} \label{zagier}
\begin{aligned}
v(1,c,m) & = (c \theta^{q^m} \tau^m -c \theta \tau^m - c^q
\tau^{m+1}, 0,\dots, 0, c\tau^{m+2} - ca_1^{q^m} \tau^{m+1}), \\
v(2,c,m) & = (c\tau^{m+1}, c\theta^{q^m}\tau^m - c\theta \tau^m -
c^q \tau^{m+1}, 0, \dots, 0, -ca_2^{q^m}\tau^{m+1}), \\
v(z,c,m) & = (0, \dots, c\tau^{m+1}, c\theta^{q^m}\tau^m - c\theta
\tau^m - c^q \tau^{m+1}, 0, \dots, 0, -ca_z^{q^m}\tau^{m+1}), \\
v(r-1,c,m) & = (0, \dots, 0, c\tau^{m+1}, c\theta^{q^m}\tau^m -
c\theta \tau^m - c^q \tau^{m+1}- ca_{r-1}^{q^m}\tau^{m+1}),
\end{aligned}
\end{equation}
where for $2 \le z \le r-2$ the possible non-zero coordinates of
$v(z,c,m)$ are $v_{z-1}$, $v_z$, and $v_{r-1}$.  Every inner
biderivation arises from an additive combination of basic $\delta_v$.

Consider $G \subseteq \Mat_{1 \times ({r-1})}(K)\{\tau\}$ consisting
of elements $u \assign (u_i)$ with $\tau$-degrees of $u_1, \dots,
u_{r-2}$ zero and the $\tau$-degree of $u_{r-1}$ less than or equal
one.  In other words, $u_1, \dots, u_{r-2} \in K$ and $u_{r-1} = c+
d\tau$ with $c, d \in K$.  Elements of $G$ give rise via
\eqref{E:MatDer} to biderivations which we will call reduced.

\begin{lemma} \label{L:GWB}
  The map $u \mapsto \delta_u \colon G \to \Der(\Psi, C)$
  induces an isomorphism
\[
G \cong \frac{\Der(\Psi, C)}{\Der_{in}(\Psi, C)} \cong
\Ext^1(E^{\vee}, C)
\]
of $\FF_q$-vector spaces.
\end{lemma}

\begin{proof}[Proof of Lemma]
  We need only to prove the first isomorphism by Lemma~\ref{L:ExtDer}.
  Let $X$ be an extension.  We want to subtract appropriate $v$'s
  in \eqref{zagier} from $\delta_X(t)= u = (u_1, u_2, \dots,
  u_{r-1})$ so that the resulting biderivation is reduced.  In this
  process, we need to keep track of the $\tau$-degrees of the $u_j$'s.
  We define the $\tau$-degree of $u$ to be the vector $d_u \assign
  (d_1(u), \dots, d_{r-1}(u))$ with $d_j(u)= \tau$-degree of $u_j$.
  Given two vectors $d$ and $d'$ with integer coefficients, we shall
  say $d \le d'$ if $d_i \le d'_i$ for all $i$.  So our claim is that
  $u$ can be reduced to a biderivation $\tilde{u}$ such that
  $d_{\tilde{u}} \le (0,\dots, 0,1)$.
  
  Let $n$ be the maximum of the integers $d_j (u)$; one has $d_j(u)\le
  n$ for all $j$.  We can modify $u$ by $v(2,c,d_1 -1)$ for an
  appropriate $c\in K$ to obtain $u'$ such that $d_1(u') < d_1(u)$ and
  $d_g(u') \le n$ for $g \ge 2$.  Subtracting an appropriate
  $v(3,c,d_2(u'))$, we obtain $u''$ such that $d_2(u'') < d_2(u')$ and
  $d_g(u'') \le n$ for $g \ge 3$.  Repeating this for $z = 4, \dots,
  r-2$ using appropriate $v(z, c, m)$, we obtain a $w \in \Mat_{1
    \times (r-1)}(K)\{ \tau \}$ whose degree vector is less than or
  equal to $(n-1, n-1,\dots, n-1, n)$.
  
  If $n\le 1$, we are done. If not ($n \ge 2$), we can subtract an
  appropriate $v(1,c,n-2)$ to obtain a vector $w$ whose degree vector
  is less than or equal to $(n-1, n-1,\dots, n-1, n-1)$.  We repeat
  the procedure in the two paragraphs above until we arrive at a
  vector $w$ whose degree vector is less than or equal to $(0,\dots,
  0, 1)$.
\end{proof}

We can now determine the $t$-module structure on $\Ext^1(E^{\vee},
C)$. By the previous lemma, it suffices to see this structure on $G$.
Consider the elements $e_i \in G \subseteq \Mat_{1 \times
  ({r-1})}(K)\{\tau\}$, $1 \le i \le r$, defined as follows.  For $1
\le i \le r-1$, we take $e_i$ to be the vector with $1$ in the $i$-th
coordinate and zeros elsewhere.  We take $e_r$ to be the vector with
$\tau$ in the last coordinate and zeros elsewhere.  The structure of a
$\FF_q[t]$-module on $G$ is completely described by the action of $t$
on elements of the form $b e_i$, $b \in K$, since additive
combinations of such elements give all of $G$.

Consider $t\ast(be_i) = (\theta + \tau) b e_i = \theta e_i + b^q \tau
e_i$.  The last is no longer an element of $G$, and we need its
equivalent vector in $G$.  Using the reduction procedure of
Lemma~\ref{L:GWB}, it is easily computed that
\[
\begin{aligned}
(\theta + \tau) (b e_{i}) & = \theta b e_n + (\textstyle{\sum_{f=0}^{r-i
   -1}} b^{q^{r-i-f}} a_{r-f}) e_r, \quad 1 \leq i \leq r-1,\\
(\theta + \tau) (b e_{r-1}) & = \theta b e_{r-1} + b^q e_r\\
(\theta + \tau) (b e_{r-2}) & = \theta b e_{r-2} + (b^{q^2} + b^q
a_{r-1}) e_r\\
(\theta + \tau) (b e_{r}) & = (\theta + \textstyle{\sum_{f=1}^{r}} b^{q^f}a_f) e_r.
\end{aligned}
\]
Thus there is a $t$-module structure $\Xi \colon \FF_q[t] \to \Mat_r(K)\{
\tau\}$ on $\Ext^1(E^{\vee}, C)$, which is completely described by
\[
  \Xi(t) \assign \begin{pmatrix}
 \theta & 0 & \cdots & \cdots & 0\\
 0  & \theta & \ddots & & 0\\
 \vdots & \ddots & \ddots & \ddots &  \vdots \\
0 & & \ddots & \theta & 0\\
\alpha_1 & \alpha_2 & \cdots & \alpha_{r-1} & \alpha_r
  \end{pmatrix}.
\]
Here $\alpha_r = \theta + a_1 \tau + a_2 \tau^2 + \cdots +
a_{r-1}\tau^{r-1} + \tau^r$, and the others are given by $\alpha_n =
\sum_{f=0}^{r-n -1} {\tau}^{r-n-f} a_{r-f}$, for $1 \le n \le r-1$.
Moreover, $\Ext^1(E^{\vee},C)$ is an extension of $\Ga^{r-1}$ by $E$.
By Lemma~\ref{L:GWB}, $(E^{\vee})^{\vee}$ is $1$-dimensional, which
completes the proof.
\end{proof}

\begin{proof}[Proof of Theorem~\ref{T:WB}(c)]
  If $F = (\Ga,\Psi)$ is a Drinfeld module, then a morphism
  $\beta\colon E \to F$ is represented by $\beta = u_0 + \dots + u_d
  \tau^d \in K\{\tau\}$ such that $\beta\Phi(a) = \Psi(a)\beta$ for
  all $a \in \FF_q[t]$.  Then $\beta$ induces an $\FF_q[t]$-module
  homomorphism $\beta^{\vee}\colon \Ext^1(F,C) \to \Ext^1(E,C)$, and
  on the level of biderivations,
\[
  (\beta^{\vee}(\delta))(a) = \delta(a) \beta, \quad \forall a \in
  \FF_q[t].
\]
Also $\beta^{\vee}$ takes $\Ext^1_0(F,C) = F^{\vee}$ into
$\Ext^1_0(E,C) = E^{\vee}$.  Since Drinfeld modules of different ranks
have no non-zero morphisms between them, we can assume that the rank
of $E$ is the same as the rank of $F$.

We continue with the considerations of the proof of part (a).  As in
\eqref{bson}, we need to measure the effect of $\beta^{\vee}$ on
biderivations in $\Der(\Psi,C)$ represented by $\alpha = b_i \tau^i$
in \eqref{E:Xred}.  We see that for $0 \leq i \leq r-1$,
\[
  \beta^{\vee}(\delta_{b_i \tau^i})(t) = b_i\tau^i\beta
    = b_i u_0^{q^i} \tau^i + \dots + b_i u_d^{q^i}\tau^{i+d}.
\]
Using \eqref{bonn} we can subtract suitable inner biderivations in
$\Der_{in}(\Phi,C)$ and find that in $\Ext^1(E,C)$,
\[
  \beta^{\vee}(\delta_{b_i \tau^i})(t) = c_{i,0}(b_i) \tau^0 + \dots +
  c_{i,r-1}(b_i)\tau^{r-1},
\]
where each $c_{i,j}(x)$ is an $\FF_q$-linear polynomial in $K[x]$
whose coefficients depend only on $\beta$.  Thus $\beta^{\vee} \colon
\Ext^1(F,C) \to \Ext^1(E,C)$ is represented by a matrix in $\Mat_{r
  \times r}(K\{\tau\})$, and $\beta^{\vee}$ restricts to a $t$-module
morphism $\beta^{\vee} \colon F^{\vee} \to E^{\vee}$.
\end{proof}

\section{Extensions of tensor powers of the Carlitz module}
\label{S:Ctens}

In this section we prove Theorem~\ref{T:Ctens}.  Recall that the
$n$-th tensor power of the Carlitz module is the $n$-dimensional pure
$t$-module $\Cn \colon \FF_q[t] \to \Mat_n(K)\{ \tau\}$ defined by
\[
  \Cn(t) \assign 
  \begin{pmatrix}
    \theta & 1 & 0 \\
    \vdots  & \ddots & 1 \\
    \tau & \cdots & \theta
  \end{pmatrix}.
\]
That is, $\Cn(t) = \theta I_n + N_n + E_n \tau$, where
\[
  N_n \assign \begin{pmatrix}
    0 & 1 & 0 \\
    \vdots & \ddots & 1 \\
    0 & \cdots & 0
    \end{pmatrix}
  \quad \textnormal{and} \quad
  E_n \assign \begin{pmatrix}
    0 & \cdots & 0 \\
    \vdots & & \vdots \\
    1 & \cdots & 0
    \end{pmatrix}.
\]
See Anderson and D.~S.~Thakur \cite{andthk90} for more details.

Fix $m < n$.  The following lemma determines representatives for
elements of $\Der(\Cm,\Cn)/\Der_{in}(\Cm,\Cn)$.  For $V \in \Mat_{n
  \times m}(K\{\tau\})$, recall the definition of $\delta_V$ from
\eqref{E:MatDer}.

\begin{lemma} \label{L:CExtgen}
Let
\[
G \assign \left\{ \begin{pmatrix} \ast &
  0 & \cdots & 0 \\ \vdots & \vdots & & \vdots \\ \ast & 0 & \cdots & 0
      \end{pmatrix} \in \Mat_{n \times m}(K) \right\}.
\]
Then $V \mapsto \delta_V \colon G \iso
\Der(\Cm,\Cn)/\Der_{in}(\Cm,\Cn)$ is an isomorphism of $\FF_q$-vector
spaces.
\end{lemma}

\begin{proof}
  Let $Q_{ij}$ be the $n \times m$ matrix with a $1$ in the $ij$-th
  entry and zeros elsewhere.  For $c \in K$ and $k \geq 0$, by taking
  $c Q_{ij} \tau^k$ for $U$ in $\delta^{(U)}(t)$, we define
\begin{multline} \label{E:Qin}
  \delta_{ijk}^{[c]} \assign \delta^{(c Q_{ij} \tau^k)}(t) =
  (cQ_{ij}\tau^k) \Cm(t) - \Cn(t)
  (cQ_{ij}\tau^k) \\
  = c(\theta^{q^k} - \theta)Q_{ij}\tau^k + c(Q_{ij}N_m -
  N_nQ_{ij})\tau^k + (c Q_{ij}E_m - c^q E_n Q_{ij})\tau^{k+1}.
\end{multline}
Since every $U \in \Mat_{n \times m}(K\{\tau\})$ is an $\FF_q$-linear
combination of matrices of the form $c Q_{ij} \tau^k$, biderivations
arising from \eqref{E:Qin} generate $\Der_{in}(\Cm,\Cn)$ as an
$\FF_q$-vector space.

Suppose that $V = (v_{ij}) \in \Mat_{n \times m}(K\{ \tau\})$ is
arbitrary and that $\deg_{\tau}(v_{ij}) \leq r$ for some $r \geq 1$.
We will show that by subtracting matrices in \eqref{E:Qin} from $V$ in
various ways we can replace $V$ by a matrix $V'$ which has each entry
of $\tau$-degree $\leq r-1$ and also $\delta_{V'}$ equivalent to
$\delta_V$ modulo $\Der_{in}(\Cm,\Cn)$.

We bootstrap our way through the entries of $V$ in the following way.
Let
\[
  D_{ij} \assign \{ \delta_{ijk}^{[c]} : k \geq 0, c \in K \}.
\]
Define a function $F$ from the set of subsets of $I \assign \{ (i,j) :
1 \leq i \leq n, 1 \leq j \leq m \}$ to itself.  We set $F(S)$ to be
those entries of $V$ whose degrees in $\tau$ can be decreased by
subtracting an element of $\sum_{(i,j) \in S} D_{ij}$, without
increasing the degrees of the other entries.  Our claim then is that
$F(I) = I$.  The following containments can be easily checked:
\begin{equation} \label{E:Fcont}
\begin{array}{lll}
F(\{2 \leq i \leq n; j=m\}) &\supseteq &\{ 2 \leq i \leq n; j=1 \}, \\
F(\{2 \leq i \leq n; j=m\} \cup \{i - j \geq 2\}) &\supseteq &\{ i - j \geq
1 \}, \\
F(\{j=m\} \cup \{i - j \geq 2\}) & \supseteq &\{ i - j
\geq 1 \} \cup \{(1,1),(1,m)\}, \\
F(\{j=m\} \cup \{i - j \geq 1\}) & \supseteq &\{ i - j
\geq 0\} \cup \{(1,m)\}, \\
F(\{j=m\} \cup \{i - j \geq -\ell \}) & \supseteq &\{ i - j
\geq -\ell - 1 \} \cup \{(1,m)\},
\end{array}
\end{equation}
where the last containment holds for all $\ell = 0, \dots, m-1$.

Therefore, we can assume that every entry of $V$ is a constant from
$K$.  Now that the $\tau$-degree of each entry of $V$ is $0$, one
checks the containments in \eqref{E:Fcont} still hold with the
exception that all sets on the right-hand side must have $\{j = 1\}$
removed.  That is, we can adjust $V$ so that it can be replaced by a
matrix in $G$, but elements of $G$ can be reduced no further.
\end{proof}

\begin{proof}[Proof of Theorem~\ref{T:Ctens}]
  Let $G$ be given as in Lemma~\ref{L:CExtgen} so that by
  Lemma~\ref{L:ExtDer}
\[
G \cong \Ext^1(\Cm,\Cn).
\]
Recalling the definition of $Q_{ij}$ from the proof above, let
\[
  e_1 \assign Q_{11}, \dots, e_n \assign Q_{n1},
\]
be basis vectors for $G$ over $K$.  Combining Lemmas~\ref{L:ExtDer}
and~\ref{L:CExtgen}, $G$ has a natural $\FF_q[t]$-module structure
which we now make explicit.  For $c \in K$,
\begin{equation} \label{E:tce1}
  t \ast (ce_1) = \Cn(t)(ce_1) = \theta c e_1 + c^q \tau e_n,
\end{equation}
and
\begin{equation} \label{E:tcei}
t \ast (ce_i) = \Cn(t)(ce_i) = \theta c e_i + c e_{i-1}, \quad 2 \leq
i \leq n.
\end{equation}
We note that as defined in \eqref{E:Qin},
\[
  \sum_{i=0}^{m-1} \delta_{n-i,m-i,0}^{[c^q]} = c^q Q_{n1} \tau -
  c^q Q_{n-m,1} = c^q\tau e_n - c^q e_{n-m},
\]
and since this sum defines an inner $(\Cm,\Cn)$-biderivation, we
subtract it from \eqref{E:tce1} and find that
\[
  t \ast (ce_1) = \theta c e_1 + c^q e_{n-m}.
\]
Therefore, combining this with \eqref{E:tcei}, we see that the
$\FF_q[t]$-module structure on $G$ can be expressed as a $t$-module by
the map $\Pi \colon \FF_q[t] \to \Mat_n(K)\{ \tau \}$ defined by
\[
  \Pi(t) \assign \left(
  \begin{array}{c|c}
  C^{\otimes(n-m)}(t) & \begin{array}{ccc}
    0 & \cdots & 0 \\ \vdots & & \vdots \\ 1 & \cdots & 0
  \end{array} 
\\ \hline
  & \\ 0 & I_m\theta + N_m \\ & \\
\end{array} \right).
\]
Thus there is an exact sequence of $t$-modules,
\[
  0 \to C^{\otimes(n-m)} \to G \to L \to 0,
\]
where $L$ is an $m$-dimensional iterated extension of $\Ga$, and one
checks that moreover $C^{\otimes (n-m)} = \Ext^1_0(\Cm,\Cn)$.
\end{proof}

\begin{remark} \label{R:nleqm}
  The $\FF_q[t]$-module structure on $\Ext^1(\Cm,\Cn)$ when $n \leq m$
  provides a different picture.  We consider the following examples.
  
  Using similar methods as in the proof of Theorem~\ref{T:Ctens} it is
  possible to show that $\Ext^1(C^{\otimes 2},C) \cong K^2$ as
  $\FF_q$-vector spaces and that the $\FF_q[t]$-module structure on
  $\Ext^1(C^{\otimes 2},C)$ is given by
\[
  t \ast \begin{pmatrix} a \\ b \end{pmatrix} =
    \begin{pmatrix} \theta & 0 \\ 1 & \theta + \tau^{-1}
    \end{pmatrix} \begin{pmatrix} a \\ b \end{pmatrix}, \quad a, b \in
    K.
\]
Thus $\Ext^1(C^{\otimes 2},C)$ is an extension of $\Ga$ by the
adjoint of the Carlitz module (see Goss \cite[\S 3.6]{goss:FF}), but
not a $t$-module.

As for an example with $n=m$, the first case to consider is
$\Ext^1(C,C)$, which itself is quite subtle.  For an inner
biderivation $\delta^{(u)} \in \Der_{in}(C,C)$ arising from $u =
c\tau^k \in K\{\tau\}$,
\[
  \delta^{(u)}(t) = c\tau^k C(t) - C(t) c\tau^k = c(\theta^{q^k} -
  \theta)\tau^k + (c - c^q)\tau^{k+1}.
\]
Unless $K$ is algebraically closed as well as perfect, it is not
possible to systematically decrease the degree in $\tau$ of
$\delta(t)$ for an arbitrary $\delta \in \Der(C,C)$ (one needs to be
able to solve equations of the form $c - c^q = \alpha$ for $\alpha \in
K$).  If $K$ is algebraically closed, the situation improves, and it
is possible to show that then $\Ext^1(C,C) \cong K$ and that in fact
$\Ext^1(C,C) \cong \Ga$ as $\FF_q[t]$-modules.
\end{remark}

\section{Periods of $t$-modules and extensions} \label{S:PerExt}

Here we would like to motivate the expectation that $\Ext^1(E,F)$, for
certain $t$-modules $E$ and $F$, can be given the structure of a
$t$-module.  Overall we have made the following observations.  If $E =
(\Ga^d,\Phi)$ is a pure $t$-module of rank $r$, then the weight of $E$
is defined to be
\[
  \wt(E) \assign d/r.
\]
What Theorems~\ref{T:WB} and~\ref{T:Ctens} have in common is that, for
certain pure $t$-modules $E$ and $F$, the $\FF_q[t]$-module
$\Ext^1(E,F)$ has the structure of a $t$-module provided $\wt(E) <
\wt(F)$ and that the submodule $\Ext^1_0(E,F)$ is itself pure and
uniformizable.

Also in the situations of Theorems~\ref{T:WB} and~\ref{T:Ctens},
$\Ext^1_0(E,F)$ behaves well under base-extension.  That is, if $L
\supseteq K$ is a perfect field, then $\Ext^1_{0/K}(E,F) \cong
\Ext^1_{0/L}(E,F)(K)$, where here $\Ext^1_{0/K}(E,F)$ is the usual
$\Ext^1_0(E,F)$ over $K$ and $\Ext^1_{0/L}(E,F)(K)$ is the set of
$K$-valued points on $\Ext^1_{0/L}(E,F)$.

These results, along with other experimental evidence, suggest that in
general $\Ext^1_0(E,F)$ is pure, uniformizable, and functorial in $K$,
for arbitrary pure uniformizable $t$-modules $E$ and $F$ with $\wt(E)
< \wt(F)$.  However, the examples in Remark~\ref{R:nleqm} show that
the situation when $\wt(E) \geq \wt(F)$ is somewhat different.

In this section we investigate the structure of $\Ext^1(E,F)$ from an
analytic point of view so as to support the veracity of the claims
above.  Our main tool will be generalizations of quasi-periodic
functions defined in \cite{bp:Gmm} and \cite{gek89} for extensions by
$\Ga$.  Moreover we study an analogue of the de Rham map of Gekeler
\cite{gek89}.  For more details on the general analytic theory of
$t$-modules see Goss \cite[Ch.~5]{goss:FF}.

Let $\KK$ be the completion of the algebraic closure of the Laurent
series field $\laurent{\FF_q}{1/\theta}$, where $\theta$ is an
independent variable.  Let $\iota\colon \FF_q[t] \to \KK$ be defined
by $t \mapsto \theta$.

Let $E = (\Ga^d, \Phi)$ and $F = (\Ga^e, \Psi)$ be uniformizable
$t$-modules over $\KK$.  We take $\Exp_E \colon \Lie(E)(\KK) \to
E(\KK)$ to be the exponential map of $E$.  Because $E$ is
uniformizable, there is an exact sequence of $\FF_q[t]$-modules
\[
  0 \to \Lambda_\Phi \to \KK^d \overset{\Exp_{\Phi}}{\to} \KK^d
  \to 0,
\]
where $\FF_q[t]$ operates by $d\Phi$ on the central $\KK^d$ and by
$\Phi$ on the right one.  Also, $\Lambda_\Phi$ is the period lattice
of $\Phi$ and is a discrete $\FF_q[t]$-submodule of $\KK^d$.
Similarly, we define $\Exp_\Psi$ and $\Lambda_\Psi$ for $F$.

\begin{lemma} \label{L:Fdel}
  For each $\delta \in \Der_0(\Phi,\Psi)$ there is a unique
  $\FF_q$-linear entire function $F_\delta \colon \KK^d \to \KK^e$
  such that for $z \in \KK^d$,
\begin{gather*}
  F_\delta(z) \equiv 0 \pmod{\deg q},\\
  F_\delta(d\Phi(a)z) = \Psi(a)F_\delta(z) + \delta(a)\Exp_\Phi(z),
  \quad \forall a \in \FF_q[t].
\end{gather*}
\end{lemma}

\begin{proof}
  The proof here is essentially the same as the that of the existence
  of the exponential function, and in particular we can easily adapt
  the proof of Proposition 2.1.4 in Anderson \cite{and86} to this
  situation.
\end{proof}

\begin{lemma} \label{L:Fin}
  Suppose $\delta^{(U)} \in \Der_0(\Phi,\Psi) \cap
  \Der_{in}(\Phi,\Psi)$.  Then
\[
  F_\delta(z) = U \Exp_\Phi(z) - \Exp_\Psi (dU \cdot z).
\]
\end{lemma}

\begin{proof}
  This follows directly from Lemma~\ref{L:Fdel}, using the fact that
  $\delta^{(U)} \in \Der_0(\Phi,\Psi)$ if and only if $dU\,d\Phi(a) =
  d\Psi(a)\,dU$ for all $a \in \FF_q[t]$.
\end{proof}

If $X = (\Ga^{d+e}, \Upsilon)$ is an extension of $E$ by $F$ defined
by a biderivation $\delta \in \Der_0 (\Phi,\Psi)$, then
Lemma~\ref{L:Fdel} and the uniqueness of the exponential function
imply that
\[
  \Exp_\Upsilon \begin{pmatrix} z \\ u \end{pmatrix} =
  \begin{pmatrix}
    \Exp_\Phi(z) \\ \Exp_\Psi(u) + F_\delta(z)
  \end{pmatrix}
\]
is the exponential function for $X$.  Also, since $E$ and $F$ are
both uniformizable, so is $X$.  Moreover, the period lattice of $X$ is
\[
  \Lambda_\Upsilon = \left\{ \begin{pmatrix} \lambda \\ \eta
  \end{pmatrix} :
  \lambda \in \Lambda_\Phi, \Exp_\Psi(\eta) + F_\delta(\lambda) = 0
  \right\}.
\]
For each $(\lambda,\eta) \in \Lambda_\Upsilon$, it follows that
$(\lambda, \eta + d\Psi(a)\mu) \in \Lambda_\Upsilon$ for all $a \in
\FF_q[t]$ and $\mu \in \Lambda_\Psi$.

\begin{proposition}
The following map is a homomorphism of $\FF_q[t]$-modules:
\[
\begin{array}{rccc}
  \DR \assign \DR_{(\Phi,\Psi)} :&
  \Der_0(\Phi,\Psi)/\Der_{si}(\Phi,\Psi) &\to &
  \Hom_{\FF_q[t]}(\Lambda_\Phi, \Psi(\KK)), \\[6pt]
  & \delta &\mapsto & (\lambda \mapsto F_\delta(\lambda)).
\end{array}
\]
\end{proposition}

\begin{proof}
  The map $\DR_0 \colon \delta \mapsto (\lambda \mapsto
  F_\delta(\lambda))$ on $\Der_0(\Phi,\Psi) \to
  \Hom_{\FF_q[t]}(\Lambda_\Phi,\Psi(\KK))$ is clearly well-defined and
  $\FF_q$-linear by Lemma~\ref{L:Fdel}.  If $\delta^{(U)} \in
  \Der_{si}(\Phi,\Psi)$, then by definition $dU = 0$, and so by
  Lemma~\ref{L:Fin}, $\DR_0(\delta^{(U)}) = 0$.  Therefore $\DR$ is
  well-defined.  Furthermore, from \eqref{E:tDer} and
  Lemma~\ref{L:Fdel}, it follows that
\[
  F_{a\ast\delta}(z) = F_{\Psi(a)\delta(\cdot)}(z) =
  \Psi(a)F_\delta(z), \quad \forall a \in \FF_q[t],
\]
and so $\DR(a\ast\delta)(\lambda) = \Psi(a)(\DR(\delta)(\lambda))$, for
all $a \in \FF_q[t]$.
\end{proof}

\begin{remark}
  In the case of extensions of Drinfeld modules by $\Ga$, the map
  $\DR$ specializes to the de Rham homomorphism of Gekeler
  \cite{gek89}.  Gekeler shows in this case that the de Rham map is in
  fact an isomorphism.  In general, determining the kernel and
  cokernel of $\DR$ is a delicate matter and will require further
  study.  However, Theorems~\ref{T:WB} and~\ref{T:Ctens} can be used
  to imply that the de Rham maps in their respective situations are
  indeed isomorphisms.  A straightforward modification of the proof of
  \cite[Thm.~3.1]{gek89} yields the following partial result.
\end{remark}

\begin{proposition} \label{P:DRinj}
  Suppose $E = (\Ga,\Phi)$ and $F = (\Ga,\Psi)$ are Drinfeld modules
  of ranks $r$ and $s$ respectively.  If $r > s$, then
  $\DR_{(\Phi,\Psi)}$ is injective.
\end{proposition}

\begin{remark} \label{R:genA}
If $A$ is the ring of functions on a smooth curve $X/\FF_q$ which are
regular away from a fixed point $\infty$, then one can also consider
Drinfeld $A$-modules (see Goss \cite[Ch.~4]{goss:FF}) and their groups
of extensions.  In \cite{gek89}, Gekeler considers $\Ext^1(E,\Ga)$,
where $E$ is a Drinfeld $A$-module of rank $r$ and shows that
\begin{equation} \label{E:ExtEGa}
  \Ext^1_0(E,\Ga) \cong \Ga^{r-1}
\end{equation}
as $A$-modules.  To prove this, Gekeler uses the de Rham isomorphism
to show that these two spaces match up exactly.  In the special case
that $A = \FF_q[t]$, the isomorphism in \eqref{E:ExtEGa} can be
generalized and proven without the de Rham isomorphism (see \cite[\S
3]{bp:Gmm}).  This is possible essentially because all $t$-modules and
biderivations of $t$-modules are determined by values of homomorphisms
on $t$ alone.

For Drinfeld $A$-modules $E$ and $F$ over a general ring $A$, one can
define biderivations and inner biderivations just as in
\S\ref{S:ExtBid}.  However, the identifications in \eqref{E:MatDer}
and thus in Lemma~\ref{L:ExtDer} are more complicated, and this
difficulty makes the $A$-module structure of $\Ext^1(E,F)$ hard to
characterize.

Although one would want to find generalizations of Theorems~\ref{T:WB}
and~\ref{T:Ctens} in the setting of Drinfeld $A$-modules, it is not
immediately clear what the precise form these analogues would take.
Since there are in general several choices of non-isomorphic rank~$1$
Drinfeld $A$-modules to consider, the exact structure of extensions by
these modules is a direction for future investigations.
\end{remark}

For the remainder of this section we will consider the implications of
the assumption that $\DR$ is an isomorphism.  In this way we attempt
to motivate the idea that $\Ext^1_0(E,F)$ can be given the structure
of a $t$-module.

Suppose $\KK^d$ is an $\FF_q[t]$-module defined by an $\FF_q$-linear
representation $\phi \colon \FF_q[t] \to \Mat_d(\KK)$ such that
$\phi(t) = \theta I_d + N$ with $N$ nilpotent.  If $\Lambda$ is any
finitely generated discrete $\FF_q[t]$-submodule of $\KK^d$ of rank
$r$, then following the language of Anderson \cite[\S 4.4]{and86}, we
call $\KK^d/\Lambda$ a \emph{$t$-torus} of dimension $d$ and rank $r$.
Through its exponential function, every uniformizable $t$-module is
isomorphic to a $t$-torus.

Let $\lambda_1, \dots, \lambda_r$ be an $\FF_q[t]$-basis for
$\Lambda_\Phi$.  Choosing this basis fixes isomorphisms
\begin{equation} \label{E:Psir}
\Hom_{\FF_q[t]}(\Lambda_{\Phi},\Psi(\KK)) \cong \Psi(\KK)^r \cong
(\KK^e)^r / (\Lambda_{\Psi})^r.
\end{equation}
We let
\begin{equation} \label{E:W}
  W \assign W_{(\Phi,\Psi)} \assign \left\{ (U\lambda_1, \dots, U\lambda_r) \in
  \Mat_{e \times r}(\KK) : \begin{array}{c}U \in \Mat_{e \times
  d}(\KK), \\ Ud\Phi(t) = d\Psi(t)U \end{array} \right\}. 
\end{equation}
Suppose $\delta^{(U)} \in \Der_0(\Phi,\Psi) \cap
\Der_{in}(\Phi,\Psi)$ for $U \in \Mat_{e \times d}(\KK)$.  By
Lemma~\ref{L:Fin} and \eqref{E:Psir},
\[
  \DR(\delta^{(U)}) = (\Exp_{\Psi}(U \lambda_1), \dots, \Exp_{\Psi}(U
  \lambda_r)) \in \Mat_{e \times r}(\KK),
\]
and so $\DR(\Der_0(\Phi,\Psi) \cap \Der_{in}(\Phi,\Psi)) =
\Exp_{\Psi}^{\oplus r}(W)$.  Thus by Corollary~\ref{C:Eseq},
$\DR_{(\Phi,\Psi)}$ induces a homomorphism of $\FF_q[t]$-modules,
\begin{equation} \label{E:Exttor}
 \DR_{(\Phi,\Psi)} : \Ext^1_0(E,F) \to \Mat_{e \times r}(\KK) / (W +
 \Lambda_{\Psi}^r),
\end{equation}
which is an isomorphism if the original $\DR_{(\Phi,\Psi)}$ is one.
Depending on $E$ and $F$, the right-hand side may or may not be
isomorphic (rigid analytically) to a $t$-torus.  We consider the
following situations.

\begin{example} \label{Ex:Drinan}
Suppose $E = (\Ga,\Phi)$ and $F = (\Ga,\Psi)$ are Drinfeld modules of
rank $2$ and $s$ respectively.  Let $\lambda_1$, $\lambda_2$ and
$\mu_1, \dots, \mu_s$ be generators for their period lattices over
$\FF_q[t]$.  The exact sequence
\[
  0 \to W \to \KK^2 \overset{f}{\to} \KK \to 0,
\]
where $f(x,y) = \lambda_1 y - \lambda_2 x$, provides a choice of
coordinates on the right-hand side of \eqref{E:Exttor}, and
$f(\Lambda_{\Psi} \oplus \Lambda_{\Psi})$ is the $\FF_q[t]$-submodule
generated by $\{\lambda_i \mu_j \}_{i = 1,2}^{j = 1, \dots, s}$.

If $s=1$, then Proposition~\ref{P:DRinj} implies that $\Ext^1_0(E,F)
\hookrightarrow \KK/(\FF_q[t]\lambda_1\mu_1 +
\FF_q[t]\lambda_2\mu_1)$, and Theorem~\ref{T:WB} (and also Woo
\cite[Prop.~7]{woo95}) confirms that this is in fact an isomorphism of
$t$-tori.  If $s \geq 2$, then it is easy to construct examples where
$f(\Lambda_\Psi \oplus \Lambda_\Psi)$ is not discrete in $\KK$, and so
in such cases the right-hand side of \eqref{E:Exttor} is not a
$t$-torus (cf.\ Theorem~\ref{L.ec}).  When we compare weights, these
observations are consistent with the discussion at the beginning of
the section, i.e.\ $\wt(E) \geq \wt(F)$ precisely when $s \geq 2$.
\end{example}

\begin{proposition} \label{P:ttorus}
  Suppose $E = (\Phi,\Ga)$ and $F = (\Psi,\Ga)$ are Drinfeld modules
  of rank $r$ and $s$ respectively, with $r > s$.  If
  $\DR_{(\Phi,\Psi)}$ is an isomorphism, then $\Ext^1_0(E,F)$ is
  isomorphic as an $\FF_q[t]$-module to a $t$-torus of dimension $r-1$
  and rank $rs$.
\end{proposition}

\begin{proof}
Let $\lambda_1, \dots, \lambda_r$ and $\mu_1, \dots, \mu_s$ be
$\FF_q[t]$-bases for $\Lambda_{\Phi}$ and $\Lambda_{\Psi}$
respectively.  Because $E$ and $F$ are both $1$-dimensional,
$\Der_{in}(\Phi,\Psi) \subseteq \Der_0(\Phi,\Psi)$.  
As in \eqref{E:W},
\[
  W = W_{(\Phi,\Psi)} =  \{ (U \lambda_1, \dots, U \lambda_r) \in
  \Mat_{1 \times r}(\KK) : U \in \KK \}. 
\]
We claim that $W \cap (\Lambda_{\Psi})^r = \{0\}$.
Suppose $w \in W \cap (\Lambda_{\Psi})^r$.  Then
\[
  w = (U\lambda_1, \dots, U\lambda_r) =  ( \nu_1, \dots, \nu_r), \quad
  \nu_i \in \Lambda_\Psi.
\]
Because $r > s$, there is a non-trivial $\FF_q[t]$-linear dependency
$\sum d\Psi(a_i)\nu_i = 0$, and thus $U\cdot\sum d\Phi(a_i)\lambda_i =
0$.  Since $\lambda_1, \dots, \lambda_r$ are linearly independent over
$\FF_q[t]$, it follows that $U = 0$.

By the following argument, the image of $\Lambda_\Psi^r$ is discrete
in $\Mat_{1 \times r}(\KK)/W$.  We proceed by induction on $r$, for
which the base case ($r=2$) is trivial.  Furthermore, it suffices to
continue with $r = s+1$, which we will now assume; the cases where
$r > s+1$ follow as straightforward consequences.  The $\KK$-linear
map
\[
\begin{array}{rccc}
  f :&
  \Mat_{1 \times r}(\KK) &\to & \Mat_{1 \times (r-1)}(\KK), \\
  & (x_1, \dots, x_r) &\mapsto & (\lambda_2 x_1 - \lambda_1 x_2,
  \dots, \lambda_r x_{r-1} - \lambda_{r-1} x_r),
\end{array}
\]
has kernel $W$.  We need to show that $V \assign \laurent{\FF_q}{1/\theta}
\cdot f(\Lambda_\Psi) \subseteq \Mat_{1 \times (r-1)}(\KK)$ is an
$rs$-dimesional vector space over $\laurent{\FF_q}{1/\theta}$.  The
typical element of $V$ has the form
\begin{equation} \label{E:nudfn}
\nu = (\mu_1, \dots, \mu_s)\, A\,\begin{pmatrix}
\lambda_2 & \cdots & 0\\
-\lambda_1 & \ddots & \vdots \\
\vdots & \ddots & \lambda_r \\
0 & \cdots & -\lambda_{r-1}
\end{pmatrix},
\end{equation}
where $A = (\alpha_{ij}) \in \Mat_{s \times
  r}(\laurent{\FF_q}{1/\theta})$.  Let $B$ denote the $r \times (r-1)$
matrix in the above formula.  By our assumption that $s = r-1$, the
matrix $AB$ is square, and its determinant is
\[
  \det(AB) =  \lambda_2 \cdots \lambda_{r-1} \sum_{j=1}^r (-1)^{r-j}
  \lambda_j \det(A_j),
\]
where $A_j$ is the $s \times s$ minor of $A$ with the $j$-th column
removed.  If $\nu = 0$, then this implies that $\det(AB) = 0$, and
since $\lambda_1, \dots, \lambda_r$ are
$\laurent{\FF_q}{1/\theta}$-linearly independent, $\det(A_j) = 0$ for each
$j$.  Thus the rank of the matrix $A$ is less than $s$, and we can
rewrite the last row of $A$ as an $\laurent{\FF_q}{1/\theta}$-linear
combination of the other rows, say $\alpha_{sj} = \sum_{i=1}^{s-1} \beta_i
\alpha_{ij}$ with $\beta_i \in \laurent{\FF_q}{1/\theta}$.  The
formulation in \eqref{E:nudfn} can be rewritten as
\[
  \nu = (\mu_1 + \beta_1 \mu_s, \dots, \mu_{s-1} + \beta_{s-1}\mu_s)
  \widetilde{A} B,
\]
where $\widetilde{A}$ is the $(s-1) \times r$ matrix obtained by
removing the last row of $A$.  Again $\nu = 0$, but then our induction
hypothesis with $r$ replaced by $r-1$ allows us to conclude that
$\widetilde{A} = 0$.
\end{proof}

\begin{remark}
  Letting $E$ and $F$ be Drinfeld modules of rank $r$ and $s$, with $r
  > s$, Proposition~\ref{P:ttorus} shows that $\Ext^1_0(E,F)$ is
  isomorphic to a $t$-torus as long as $\DR_{(\Phi,\Psi)}$ is an
  isomorphism.  Since in this case $\wt(E) < \wt(F)$, our discussion
  at the beginning of this section leads us to speculate that
  $\DR_{(\Phi,\Psi)}$ is an isomorphism and that the $t$-torus
  isomorphic to $\Ext^1_0(E,F)$ is in fact a pure uniformizable
  $t$-module.  In addition, using the techniques of
  Theorem~\ref{T:WB}, it is possible to show that $\Ext^1_0(E,F)$ is
  isomorphic to a $t$-module of dimension $r-1$, though we omit the
  details.
\end{remark}

\section{Elliptic curves} \label{licht}

Let $E_1$ and $E_2$ be elliptic curves over $\CC$. Let $G \assign
\Ext^1_{\CC}(E_1, E_2)$ be the extension group in the category of complex
abelian varieties. By the Poincar\'e reducibility theorem, this is a
torsion group.  The set of complex points $E(\CC)$ of an elliptic
curve $E$ over $\CC$ can be viewed as a complex torus written
$E^{an}$; so one may also consider $A \assign \Ext^1(E_1^{an},
E_2^{an})$, the extension group in the category of complex tori. There
is a natural homomorphism $G \to A$; the image is the torsion subgroup
of $A$ (see \cite[Rmk.~6.2, p.~23]{bl}).

We present next a theorem of Lichtenbaum (1960's, unpublished);
however, our proof is different in that it is based on periods.  By
comparison to Example~\ref{Ex:Drinan} and Proposition~\ref{P:ttorus},
the theorem below indicates that the situation for elliptic curves
runs quite parallel to our own.  Specifically, in light of in
Example~\ref{Ex:Drinan}, it suggests that $\Ext^1(E,F)$ will rarely be
representable as a $t$-module for nonisogenous Drinfeld modules $E$
and $F$ of rank 2.  Moreover, it suggests that $\Ext^1(E,F)$, for
general $t$-modules $E$ and $F$, will not always have the structure of
a $t$-module.

\begin{theorem}[Lichtenbaum] \label{L.ec}
Let $E_1$ and $E_2$ be elliptic curves over $\CC$, and let $G \assign
\Ext^1_{\CC}(E_1,E_2)$ and $A \assign \Ext^1(E_1^{an}, E_2^{an})$.
\begin{enumerate}
\item[(a)] If $E= E_1 = E_2$ has complex multiplication, then $A$ is
  naturally isomorphic to $E(\CC)$ as abelian groups.  Under this
  isomorphism, the group $G$ is identified with the torsion points of
  $E(\CC)$.
  
\item[(b)] If $E_1$ and $E_2$ are isogenous and admit complex
  multiplication, then $A$ is isogenous to $E_1(\CC)$ and $E_2(\CC)$.
  
\item[(c)] If at least one of $E_1$ and $E_2$ does not admit complex
  multiplication, then the natural topology on $G$ is non-Hausdorff.
  Therefore, $G$ is not the set of complex points of a complex
  algebraic variety with the classical topology.
\end{enumerate}
\end{theorem}

\begin{proof}
  Given $E_i$ ($i =1,2$) as the quotient of $\CC$ by lattices $\ZZ +
  \ZZ \tau_i$, where each $\tau_i$ may be taken to be have positive
  imaginary part, we see from \cite[Prop.~5.7, p.~21]{bl} that $A$ is
  naturally the quotient of $\CC$ by the subgroup $\Lambda$ generated
  by $1, \tau_1, \tau_2$ and $\tau_1\tau_2$.
  
  Let us suppose that an elliptic curve $X$ has complex
  multiplication; let us think of $X(\CC)$ as a quotient of $\CC$ by
  $\ZZ + \ZZ \tau$.  By the theory of complex multiplication
  \cite{ser1}, we have the following: (i) $\tau$ lies in an imaginary
  quadratic field; (ii) $\tau^2 = a \tau + b$ for some integers $a$
  and $b$; and (iii) if $X'$ (with $X'(\CC)= \CC/{(\ZZ + \ZZ \tau')}$) is 
an elliptic curve isogenous to $X$, then $\tau'$ and $\tau$ lie in the same
  imaginary quadratic field and $X'$ also has complex multiplication
  with the same CM-field.
  
  If $X(\CC) = E = E_1 = E_2$, then we have $\tau = \tau_1 = \tau_2$.
  By (i) and (ii), we get that $\Lambda$ is the lattice $\ZZ + \ZZ
  \tau$.  Thus we obtain that $A$ is naturally identified with
  $E(\CC)$.  The last statement in (a) follows from \cite[Rmk.~6.2,
  p.~23]{bl}.
  
  If $E_1$ and $E_2$ are isogenous curves with complex multiplication,
  then $\tau_1$ and $\tau_2$ lie in the same imaginary quadratic field
  $F$; in this case, $\Lambda$ is isomorphic to a fractional ideal of
  an order of $K$. This proves (b).
  
  For (c), suppose at least one of $E_1$ and $E_2$ does not admit
  complex multiplication.  By the fundamental theorem of complex
  multiplication, $\tau_1$ and $\tau_2$ are not both contained in one
  imaginary quadratic field.  In other words, in this case $\Lambda$
  is a subgroup of $\ZZ$-rank greater than two, and so it is not a
  discrete subgroup of $\CC$.
\end{proof}

\begin{remark} (C.~Schoen) (a) If $A$ and $B$ are abelian varieties defined 
  over $\overline{\QQ}$, then the natural map
  $\Ext^1_{\overline{\QQ}}(A,B) \to \Ext^1_{{\CC}}(A,B)$ is an
  isomorphism.
  
  (b) Any abelian surface corresponding to an element of
  $\Ext^1_{\CC}(E_1, E_2)$ is isogenous to the product $E_1 \times
  E_2$.  Any complex abelian variety, which is isogenous to a product
  of CM-elliptic curves, is itself a product of CM-elliptic curves
  \cite{sch}.
\end{remark}   

\section{Extensions of $t$-motives} \label{S:Rmk}

In this section we explore extensions from the standpoint of
$t$-motives and examine avenues for further study.  Given two
$t$-modules $E$, $F$ over $K$, we can consider the associated Anderson
$t$-motives $M(E)$ and $M(F)$ \cite{and86}.  Since the functor $M$
which sends a $t$-module to its associated $t$-motive is
contravariant, we obtain a map
\[
  M^*: \Ext^1(E,F) \to\Ext^1((M(F), M(E));
\]
because $M$ gives an anti-equivalence of categories of $t$-modules and
$t$-motives \cite[Thm.~1]{and86}, $M^*$ is an isomorphism.  If one is
interested in computing just the group of extensions of $t$-modules,
then it is relatively easy to compute extensions in the category of
$t$-motives.  We formulate this precisely in the next lemma.

The evident functor $f$ from the category $\mathcal{T}$ of $t$-motives
to the category $\mathcal{C}$ of left $K[t, \tau]$-modules is fully
faithful \cite[\S 1.2]{and86}.  Here $K[t, \tau]$ is the
noncommutative ring generated by $t$ and $\tau$ with the relations,
$t\tau = \tau t, xt = tx, \tau x = x^q\tau$, for all $x\in K$.

\begin{lemma} \label{L:tmot}
For any $t$-motives $A$ and $B$ over $K$, the natural map 
\[
  f^*: \Ext^1_{\mathcal{T}}(A,B) \to \Ext^1_{\mathcal{C}}(f(A), f(B))
\]
is an isomorphism.
\end{lemma}

\begin{proof} The fact that $f$ is fully faithful implies that
  $f^*$ is injective in the following way.  Suppose the images of
  $\alpha$ and $\beta$ under $f^*$ coincide.  Pick representatives of
  $\alpha$ and $\beta$, i.e.\ extensions $X_1$ and $X_2$ of $A$ by $B$
  which satisfy $f(X_1) \cong f(X_2)$.  We obtain a commutative
  diagram
\[
\xymatrix{ 0 \ar[r] & f(B) \ar[r] \ar@{=}[d] & f(X_1)
  \ar[r] \ar[d]_{\gamma} & f(A) \ar[r] \ar@{=}[d] & 0 \\
  0 \ar[r] & f(B) \ar[r] & f(X_2) \ar[r] & f(A) \ar[r] & 0,}
\]
where $\gamma$ is an isomorphism of $K[t, \tau]$-modules. Since $f$
is fully faithful, $\gamma$ is an isomorphism of $t$-motives. More
precisely, $\gamma$ is induced by an isomorphism $X_1 \iso X_2$ of
$t$-motives.
 
It remains to show the surjectivity of $f^*$.  For this, we have to
show the following: given any extension $X$ of $f(A)$ by $f(B)$ in
$\mathcal C$, the left $K[t, \tau]$-module $X$ is a $t$-motive, i.e.\ 
(i) it is free and finitely generated as a $K[\tau]$-module and (ii)
the associated primes of $X' \assign X/{{\tau}X}$, viewed as a module
over the commutative ring $R \assign K[t]$, consist only of the
principal ideal $I \assign (t-\theta)$.  Geometrically, we want the
coherent sheaf associated to $X'$ on $\A^1$ to be supported only at
the point $t=\theta$.

Condition (i) is clearly satisfied by $X$ by general properties of
modules over the ring $K[\tau]$ \cite[Prop.~5.4.9]{goss:FF}.  For
(ii), consider an extension $Q$ of $P$ by $N$ where $P$ and $N$ are
finitely generated $R$-modules.  Every associated prime of $Q$ is an
associated prime of either $P$ or $N$.  Also, if $H$ is a quotient
module of $N$, then every associated prime of $H$ is an associated
prime of $N$ (this is easy to see via the geometric interpretation).

Now, by assumption, the associated primes of the $R$-modules $A' \assign
f(A)/{{\tau f(A)}}$ and $B' \assign f(B)/{{\tau f(B)}}$ consist of
just the ideal $I$.  The $R$-module $X'$ is an extension of $A'$ by
$B''$ ($=$ a quotient module of $B'$). So we may apply the
comments in the previous paragraph to the extension $X'$ to deduce
that $X'$ satisfies (ii).
\end{proof}

\begin{remark}
  Lemma~\ref{L:tmot} shows that extensions may be computed via
  resolutions of $t$-motives by free $K[t, \tau]$-modules.
  Furthermore, it implies the injectivity of the map
\[
  f^*: \Ext^2_{\mathcal{T}}(A,B) \to \Ext^2_{\mathcal{C}}(f(A), f(B)).
\]
For $M$ and $N$ as in Remark~\ref{van} (and $K = \FF_q$), this gives 
$\Ext^2_{\mathcal{T}}(M,N) = 0$.
\end{remark}
  
\begin{remark}[Analogy with $\mathcal{D}$-modules]
  A.~Rosenberg points out that $K[t, \tau]$ and, more generally, skew
  polynomial algebras are analogous to Weyl algebras in $\mathcal
  D$-module theory (the first Weyl algebra is $\CC[x,\partial]$ with
  the relation $\partial x - x \partial = 1$) in that they are all
  special cases of \emph{hyperbolic} algebras \cite[Ch.~II]{ros}.
  Since $\Ext$'s of certain (but not all) $\mathcal{D}$-modules
  possess a nontrivial structure of a $\mathcal{D}$-module, one
  may expect the same to be true for $t$-modules.
  
  The analogy with $\mathcal D$-modules is best viewed within the
  context of opers (see \cite[\S 7.3.14]{bzf}; a Drinfeld module is
  an example of a Frobenius oper) and noncommutative algebraic
  geometry (see \cite[Rmk.~5.3.5, \S 6]{bstu}, \cite[\S 0.6]{car},
  \cite{ros}, \cite{tag93}). The analogy between Drinfeld modules and
  noncommutative tori is explained in \cite{man}.
\end{remark}

\begin{remark} \label{van}
  Geometrically interpreting the definition of a $t$-motive
  \cite[\S 1.2]{and86}, a $t$-motive $M$ is a special sheaf
  $\mathcal{F}_M$ over a noncommutative surface $S$ \cite{ros} given
  by the product of the noncommutative affine line $\A^1_{nc}$ by the
  commutative affine line ${\A}^1$.  We think of $S$ as fibered over
  $\mathbb{A}^1$ (viewed horizontally) with (vertical) fibers
  $\mathbb{A}^1_{nc}$.  The sheaf $\mathcal{F}_M$ is the ideal sheaf
  of a curve $X_M \subseteq S$ for which the following hold.
  
  (a) $X_M$ is finite over the two axes, via the projections to the
  components; in other words, $X_M$ is transversal to the horizontal
  and vertical fibers.
  
  (b) The intersection of $X_M$ with the horizontal axis (which
  corresponds to $\tau =0$) is a nilpotent subscheme of $\mathbb{A}^1$
  supported at the point $t=\theta$; (a) assures us that the
  intersection is a proper subscheme of $\mathbb{A}^1$.
  
  If $M$ and $N$ are distinct $t$-motives, then the group
  $\Ext^1_{\mathcal{T}}(M,N)$ can be interpreted via the
  ``intersection scheme'' of the curves $X_M$ and $X_N$ in $S$. We can
  view the subgroup $\Ext^1_{0}(M,N)$ (cf.\ Corollary~\ref{C:Eseq}) of
  $\Ext^1_{\mathcal{T}}(M,N)$ as the nontrivial part of the
  intersection locus corresponding to points distinct from $t= \theta,
  \tau = 0$, as in (b).
  
  The geometric situation is especially clear in the case $K = \FF_q$:
  $S$ is the usual commutative affine plane $\A^2$; if $g(t, \tau)$
  and $h(t, \tau)$ (assumed to have no common factors) are defining
  equations for the curves $X_M$ and $X_N$, then
  $\Ext^1_{\mathcal{T}}(M,N)$ is isomorphic to the quotient module
  $\FF_q[t,\tau]/{(g,h)}$ of the commutative ring $\FF_q[t,\tau]$.
\end{remark}

\begin{remark}
  The results and ideas in this paper are used in an ongoing project
  with Thakur, whose aim is to relate extension groups of $t$-modules
  to values of zeta functions in the spirit of \cite{andthk90},
  \cite{del2}, \cite{mr}. 
\end{remark}

\section*{Acknowledgements}

We wish to thank Greg Anderson, Dale Brownawell, David Goss, Steve
Lichtenbaum, and Dinesh Thakur for encouragement, useful discussions,
and helpful advice.  We also thank the referee for many helpful
comments and suggestions.  Part of the work on this paper was done
when the second author was visiting MPIM (Bonn), which he thanks for
its hospitality, support, and excellent working conditions.

\end{document}